\documentclass[amssymb,12pt]{amsart}
\usepackage{latexsym}
\usepackage{amsfonts}
\usepackage{amsmath}
\usepackage{amssymb}
\usepackage{mathtext}
\usepackage[mathscr]{eucal}
\newtheorem{propo}{Proposition}

\newtheorem{lemma}[propo]{Lemma}

\newtheorem{theor}[propo]{Theorem}

\newtheorem{prob}[propo]{Problem}

\newcommand{\Ker}{\operatorname{Ker}}

\newcommand{\Irr}{{\mathrm {Irr}}}

\newcommand{\IBRL}{{\mathrm {IBr}}_{\ell}}

\newcommand{\Ind}{{\mathrm {Ind}}}

\newcommand{\Hom}{{\mathrm {Hom}}}

\newcommand{\ZZ}{{\mathbb Z}}

\newcommand{\FF}{{\mathbb F}}

\newcommand{\EC}{\mathcal{E}}

\newcommand{\GC}{\mathcal{G}}

\newcommand{\ta}{\hspace{0.5mm}^{2}\hspace*{-0.2mm}}
\newcommand{\tb}{\hspace{0.5mm}^{3}\hspace*{-0.2mm}}
\title[Cross characteristic representations of $\tb D_{4}(q)$]
{Cross Characteristic Representations of $\tb D_{4}(q)$ are Reducible
over Proper Subgroups}

\author{}
\date{Nov. 27, 2007}

\begin{document}

\maketitle

\bigskip
\centerline{by}
\bigskip

\centerline{Hung Ngoc Nguyen and Pham Huu Tiep}
\centerline{Department of Mathematics}
 \centerline{University of Florida}
 \centerline{Gainesville, FL 32611, U.S.A.}
\centerline{hnguyen@math.ufl.edu, tiep@math.ufl.edu}
\bigskip

\bigskip
\centerline{with an appendix by}
\bigskip

\centerline{Frank Himstedt}
\centerline{Technische Universit\"at M\"unchen}
 \centerline{Zentrum Mathematik -- M11}
 \centerline{Boltzmannstr. 3, 85748 Garching, Germany}
\centerline{himstedt@ma.tum.de}

\bigskip
\bigskip
\noindent Running title: Cross characteristic representations of $\tb D_{4}(q)$
\bigskip


 \bigskip
\bigskip
\noindent
{\bf Abstract.} We prove that the restriction of any absolutely irreducible
representation of Steinberg's triality groups $^3D_4(q)$ in
characteristic coprime to $q$ to any proper subgroup is reducible.

\bigskip
\noindent Part of the paper was written while the authors were
participating in the Special Semester on Group Representation Theory
at the Bernoulli Center, Ecole Polytechnique Federale de Lausanne
(EPFL), Switzerland. It is a pleasure to thank the organizers,
Professors. M. Geck, D. Testerman, and J. Thevenaz for generous hospitality and
support and stimulating environment. The authors are grateful to J. Saxl for
helpful discussion.

\noindent The second author gratefully acknowledges the support of the NSF (grant
DMS-0600967).

\newpage

\section{Introduction}

Finite primitive permutation groups have been studied since the pioneering
work of Galois and Jordan on group theory; they have had important
applications in many different areas of mathematics.

If $G$ is a primitive permutation group with a point stabilizer $M$ then
$M<G$ is a maximal subgroup. Thanks to work of Aschbacher, O'Nan, Scott \cite{AS}, and
of Liebeck, Praeger, Saxl, and Seitz \cite{LPS}, \cite{LiS}, most problems involving
such a $G$ can be reduced to the case where $G$ is a finite classical group. In this case,
Aschbacher's theorem \cite{A} describes all possible choices for the maximal subgroup $M$ of
$G$. Work of Kleidman and Liebeck \cite{KL} and others essentially reduces the question of
whether a subgroup $M$ from this list is indeed maximal in $G$ to the following problem.

\begin{prob}\label{restr}
{\sl Let $\FF$ be an algebraically closed field of characteristic
$\ell$. Classify all triples $(K,V,H)$ where $K$ is a finite group
with $K/Z(K)$ almost simple, $V$ is an $\FF K$-module of dimension
greater than one, and $H$ is a proper subgroup of $K$ such that the
restriction $V|_{H}$ is irreducible.}
\end{prob}

Our main focus is on the case where $K$ is a finite group of Lie type in characteristic
$\neq \ell$. If, furthermore, $K$ is of type $A$, then Problem \ref{restr} has been
solved recently in \cite{KT}. On the other hand, the case where $K$ is an exceptional
group of type $G_{2}$, $\ta B_{2}$, or $\ta G_{2}$, was settled in \cite{N}.

The main result of this paper is the following:

\begin{theor}\label{main}
{\sl Let $G = \tb D_{4}(q)$ and let $\Phi$ be any irreducible
representation of $G$ in characteristic $\ell$ coprime to $q$.
If $H$ is any proper subgroup of $G$ and $\deg(\Phi) > 1$, then $\Phi|_{H}$
is reducible.}
\end{theor}

In the case of complex representations, Theorem \ref{main} was proved by
Saxl \cite{S}.

\section{Basic Reduction}
Given a finite group $X$,
we denote by $\mathfrak{d}_{\ell}(X)$ and $\mathfrak{m}_{\ell}(X)$
the smallest and largest degrees of absolutely irreducible
representations of degree greater than one of $X$ in characteristic
$\ell$; furthermore, let
$\mathfrak{m}_{\mathbb{C}}(X)=\mathfrak{m}_{0}(X)$. From
now on, $\FF$ stands for an algebraically closed field of
characteristic $\ell$, and $q$ is a power of a prime $p$.
If $\chi$ is a complex character of $X$, we
denote by $\widehat{\chi}$ the restriction of $\chi$ to
$\ell$-regular elements of $X$.
By $\IBRL(X)$ we denote the set of irreducible $\ell$-Brauer characters, or
the set of absolutely irreducible $\FF X$-representations, depending on the context.

First we record a few well-known statements.

\begin{lemma}
\label{lemma0} {\sl Let $K$ be a finite group. Suppose $V$ is an
irreducible $\FF K$-module of dimension greater than one, and $H$ is
a proper subgroup of $K$ such that the restriction $V|_H$ is
irreducible. Then
$$\sqrt{|H/Z(H)|} \geq \mathfrak{m}_{\mathbb{C}}(H)\geq\mathfrak{m}_{\ell}(H)\geq \dim
(V)\geq \mathfrak{d}_{\ell}(K).$$}
\end{lemma}

\begin{lemma} {\rm \cite[p. 190]{I}}
\label{lemmaIsaacs} {\sl  Let $K$ be a finite group and $H$ be a
normal subgroup of $K$. Let $\chi\in \Irr(K)$ and $\theta\in\Irr(H)$
be a constituent of $\chi_H$. Then $\chi(1)/\theta(1)$ divides
$|K/H|$.}
\end{lemma}

\begin{lemma}
\label{lemma1} {\sl Let $K$ be a simple group and $V$ an absolutely
irreducible $\FF K$-module of dimension greater than one. Suppose
$H$ is a subgroup of $K$ such that $V|_H$ is irreducible. Then
$Z(H)=C_K(H)=1$.}
\end{lemma}

In the following theorem, we use results and notation of \cite{K}.

\begin{theor}[Reduction Theorem]
\label{reduction} {\sl Let $G = \tb D_{4}(q)$ and let $\varphi$ be
an irreducible representation of $G$ in characteristic $\ell$
coprime to $q$. Suppose $\varphi(1)>1$ and $M$ is a maximal subgroup
of $G$ such that $\varphi|_{M}$ is irreducible. Then $M$ is
$G$-conjugate to one of the following groups:

{\rm (i)} $P$, a maximal parabolic subgroup of order $q^{12}(q^{6}-1)(q-1)$,

{\rm (ii)} $Q$, a maximal parabolic subgroup of order $q^{12}(q^{3}-1)(q^{2}-1)$,

{\rm (iii)} $G_2(q)$,

{\rm (iv)} $\tb D_4(q_0)$ with $q=q_0^2$.}
\end{theor}

\begin{proof}
By \cite[Theorem 4.1]{MMT}, $\mathfrak{d}_{\ell}(^3D_4(q))\geq q^5-q^3+q-1$ for every
$\ell$ coprime to $q$. Next, according to \cite{K}, if $M$ is a maximal subgroup of $G$, but
$M$ is not a maximal parabolic subgroup, then $M$ is
$G$-conjugate to one of the following groups:

1) $G_2(q)$,

2) $PGL^\epsilon_3(q)$, where $4\leq q\equiv\epsilon 1(\bmod$ $3)$,
$\epsilon=\pm$,

3) $^3D_4(q_0)$ with $q=q_0^\alpha$, $\alpha$ prime, $\alpha\neq3$,

4) $L_2(q^3)\times L_2(q)$, where $2\mid q$, a fundamental subgroup,

5) $C_G(s)=(SL_2(q^3)\circ SL_2(q)).2$, $q$ odd, involution
centralizer,

6) $((\mathbb{Z}_{q^2+q+1})\circ SL_3(q)).f_+.2$, where
$f_+=(3,q^2+q+1)$,

7) $((\mathbb{Z}_{q^2-q+1})\circ SU_3(q)).f_-.2$, where
$f_-=(3,q^2-q+1)$,

8) $(\mathbb{Z}_{q^2+q+1})^2.SL_2(3)$,

9) $(\mathbb{Z}_{q^2-q+1})^2.SL_2(3)$,

10) $\mathbb{Z}_{q^4-q^2+1}.4$.

We only need to consider the following cases.

$\bullet$ $M=PGL^\epsilon_3(q)$, where $4\leq q\equiv\epsilon
1(\bmod$ $3)$, $\epsilon=\pm$.

We have $|PGL^\epsilon_3(q)|=q^3(q^2-1)(q^3\pm 1)$. So
$\mathfrak{m}_\mathbb{C}(M)\leq \sqrt{q^3(q^2-1)(q^3+1)}
 < q^5-q^3+q-1$ for every $q\geq4$. Therefore
$\mathfrak{m}_\mathbb{C}(M)<\mathfrak{d}_{\ell}(^3D_4(q))$, contradicting
Lemma \ref{lemma0}.

$\bullet$ $M= \ ^3D_4(q_0)$ with $q=q_0^\alpha$, $\alpha$ prime,
$\alpha\neq2$, $3$.

We have $|^3D_4(q_0)|=q_0^{12}(q_0^6-1)^2(q_0^4-q_0^2+1)<q_0^{28}$.
Hence, $\mathfrak{m}_{\mathbb{C}}(M)\leq\sqrt{q_0^{28}}=q_0^{14}$.
Since $\alpha$ is prime and $\alpha\neq 2$, $3$, $\alpha\geq5$. It
follows that
$\mathfrak{m}_{\mathbb{C}}(M)<q^{14/5}<q^5-q^3+q-1\leq\mathfrak{d}_{\ell}(^3D_4(q))$.

$\bullet$ $M=L_2(q^3)\times L_2(q)$, where $2\mid q$.

It is well known that $\mathfrak{m}_{\mathbb{C}}(L_2(q))=q+1$ except
that $\mathfrak{m}_{\mathbb{C}}(L_2(2))=2$,
$\mathfrak{m}_{\mathbb{C}}(L_2(3))=3$ and
$\mathfrak{m}_{\mathbb{C}}(L_2(5))=5$. So we have
$\mathfrak{m}_{\mathbb{C}}(L_2(q^3))=q^3+1$ for every $q$. Hence
$\mathfrak{m}_{\mathbb{C}}(M)= (q+1)(q^3+1)$ for $q\geq 4$ and
$\mathfrak{m}_{\mathbb{C}}(M)=18$ for $q=2$. It is easy to see that
$(q+1)(q^3+1)<q^5-q^3+q-1$ for every $q\geq 4$. When $q=2$, we also
have $\mathfrak{m}_{\mathbb{C}}(M)=18<25=2^5-2^3+2-1$. Therefore,
$\mathfrak{m}_{\mathbb{C}}(M)<\mathfrak{d}_{\ell}(^3D_4(q))$ for
every $q$.

$\bullet$ $M=C_G(s)=(SL_2(q^3)\circ SL_2(q)).2$, $q$ odd, $s$ an involution.

Here $C_{G}(M) \ni s \neq 1$, contradicting Lemma \ref{lemma1}.

$\bullet$ $M=((\mathbb{Z}_{q^2+q+1})\circ SL_3(q)).f_+.2$, where
$f_+=(3,q^2+q+1)$.

By Lemma \ref{lemmaIsaacs}, $\mathfrak{m}_{\mathbb{C}}(M)\leq
2.f_+.\mathfrak{m}_{\mathbb{C}}((\mathbb{Z}_{q^2+q+1})\circ
SL_3(q))\leq 2.f_+.\mathfrak{m}_{\mathbb{C}}(SL_3(q))$. From
\cite{SF}, we have

\begin{center}
\begin{tabular}{c}
$\mathfrak{m}_{\mathbb{C}}(SL_3(q))= \left\{\begin{tabular}{ll}
$8$, & $q=2$,\\
$39$, & $q=3$,\\
$84$, & $q=4$,\\
$(q+1)(q^2+q+1)$, & $q\geq5$.\\
 \end{tabular} \right.$
\end{tabular}
\end{center}
It is easy to check that $2.f_+.\mathfrak{m}_{\mathbb{C}}(SL_3(q))<
q^5-q^3+q-1$ for every $q\geq2$. Therefore
$\mathfrak{m}_{\mathbb{C}}(M)<q^5-q^3+q-1$.

$\bullet$ $M=((\mathbb{Z}_{q^2-q+1})\circ SU_3(q)).f_-.2$, where
$f_-=(3,q^2-q+1)$.

By Lemma \ref{lemmaIsaacs}, $\mathfrak{m}_{\mathbb{C}}(M)\leq
2.f_-.\mathfrak{m}_{\mathbb{C}}((\mathbb{Z}_{q^2-q+1})\circ
SU_3(q))\leq 2.f_-.\mathfrak{m}_{\mathbb{C}}(SU_3(q))$. From
\cite{SF}, we have

\begin{center}
\begin{tabular}{c}
$\mathfrak{m}_{\mathbb{C}}(SU_3(q))= \left\{\begin {array}{ll}
8, & q=2,\\
(q+1)^2(q-1), & q\geq3.\\
 \end {array} \right.$
\end{tabular}
\end{center}
It is easy to check that $2.f_-.\mathfrak{m}_{\mathbb{C}}(SU_3(q))<
q^5-q^3+q-1$ for every $q\geq3$. Therefore
$\mathfrak{m}_{\mathbb{C}}(M)<q^5-q^3+q-1$ for every $q\geq3$.

When $q=2$, we have $\mathfrak{m}_{\mathbb{C}}(M)\leq
\sqrt{|M|}=36$. Therefore if $\varphi|_M$ is irreducible then $\deg
(\varphi)\leq 36$. Inspecting the character tables of $^3D_4(2)$ in
\cite{Atlas1} and \cite{Atlas2}, we see that
$\deg(\varphi)=25$ for $\ell=3$ or $\deg(\varphi)=26$ for
$\ell\neq3$. Moreover, when $\ell\neq3$, $\varphi$ is the reduction
modulo $\ell$ of the unique irreducible complex representation $\rho$ of
degree $26$. Since $26\nmid|M|=1296$, $M$ does not have any
irreducible complex representation of degree $26$, whence
$\rho|_{M}$ and $\varphi|_{M}$ must be
reducible. When $\ell=3$, $M=((\mathbb{Z}_{3})\circ
SU_3(2)).3.2\simeq3^{1+2}.2S_4$. So
$\mathfrak{m}_{3}(M)=\mathfrak{m}_{3}(3^{1+2}.2S_4)=\mathfrak{m}_{3}(2S_4)\leq
\mathfrak{m}_{\mathbb{C}}(2S_4)\leq\sqrt{|2S_4|}<7$. Therefore if
$\deg(\varphi)=25$ then $\varphi|_M$ is reducible.

$\bullet$ $M=(\mathbb{Z}_{q^2\pm q+1})^2.SL_2(3)$.

We have $\mathfrak{m}_{\mathbb{C}}(M)\leq |SL_2(3)|=24$. Since
$q^5-q^3+q-1>24$ for every $q\geq 2$,
$\mathfrak{m}_{\mathbb{C}}(M)<\mathfrak{d}_{\ell}(^3D_4(q))$.

$\bullet$ $M=\mathbb{Z}_{q^4-q^2+1}.4$.

We have $\mathfrak{m}_{\mathbb{C}}(M)\leq 4 < \mathfrak{d}_{\ell}(G)$. The Theorem is proved.
\end{proof}

\section{Restrictions to $G_{2}(q)$ and to $\tb D_{4}(\sqrt{q})$}
In this section we handle two of the maximal subgroups singled out in
Theorem \ref{reduction}.

\begin{theor}\label{g2}
{\sl Let $M \simeq G_{2}(q)$ be a subgroup of $G = \tb D_{4}(q)$ and $\varphi \in \IBRL(G)$
be of degree $>1$. Then $\varphi|_{M}$ is reducible.}
\end{theor}

\begin{proof}
Assume the contrary: $\varphi|_{M}$ is irreducible. By Lemma
\ref{lemma0}, $\varphi(1) < \sqrt{|M|} < q^{7}$. We will identify
the dual group $G^{*}$ with $G$. By the fundamental result of
Brou\'e and Michel \cite{BM}, $\varphi$ belongs to a union
$\EC_{\ell}(G,s)$ of $\ell$-blocks, labeled by a semisimple
$\ell'$-element $s \in G$. Moreover, by \cite{HM}, $\varphi(1)$ is
divisible by $(G:C_{G}(s))_{p'}$. Assume $s \neq 1$. Then it is easy
to check, using \cite{DM} for instance, that $(G:C_{G}(s))_{p'} \geq
q^{8}+q^{4}+1$. Since $\varphi(1) < q^{7}$, it follows that $s = 1$,
i.e. $\varphi$ belongs to a unipotent block.

According to \cite{K}, $M = C_{G}(\tau)$ for some (outer) automorphism $\tau$ of order $3$ of
$G$. Furthermore, the degrees of all complex irreducible characters of
$G$ are listed in \cite{DM}. An easy inspection reveals that $G$ has a unique irreducible
character of degree $\psi(1)$ for every unipotent character $\psi \in \Irr(G)$. It follows
that every unipotent (complex) character of $G$ is $\tau$-invariant.

Next we show that $\varphi$ is also $\tau$-invariant.
First consider the case where $q$ is odd. Then Corollary 6.9 of \cite{G} states that the
$\ell$-modular decomposition matrix of $G$ has a lower unitriangular
shape. In particular, this implies that $\varphi$ is an integral
linear combination of $\hat{\psi}$, with $\psi \in \EC(G,1)$, the
set of unipotent characters of $G$. But each such $\psi$ is
$\tau$-invariant, whence $\varphi$ is $\tau$-invariant.
Now assume that $q$ is even. Then $\ell \neq 2$, and so it is a good
prime for $\GC$, and $\ell$ does not divide $|Z(\GC)|$, where $\GC$
is the simple, simply connected algebraic group of type $D_{4}$.
Hence, by the main result of \cite{GH}, $\{\hat{\psi} \mid \psi \in
\EC(G,1)\}$ is a basic set of Brauer characters of
$\EC_{\ell}(G,1)$. It follows that $\varphi$ is an integral linear
combination of $\hat{\psi}$, with $\psi \in \EC(G,1)$, and so it is
$\tau$-invariant as above.

Consider the semidirect product $\tilde{G} = G \cdot \langle \varphi \rangle$. Then
$G \lhd \tilde{G}$, and $\tilde{G}/G$ is cyclic. Since $\varphi$ is $\tilde{G}$-invariant,
it extends to $\tilde{G}$ by \cite[Theorem III.2.14]{F}. But $C_{\tilde{G}}(M) \ni \tau \neq 1$,
hence $\varphi|_{M}$ cannot be irreducible by Lemma \ref{lemma1}.
\end{proof}

\begin{theor}\label{square}
{\sl Let $H \simeq \tb D_{4}(q)$ be a maximal subgroup of $G = \tb D_{4}(q^{2})$ and
$V \in \IBRL(G)$ be of dimension $>1$. Then $V|_{H}$ is reducible.}
\end{theor}

\begin{proof}
Again assume the contrary. We consider a long-root parabolic subgroup
$P = q^{2+16} \cdot SL_{2}(q^{6}) \cdot \ZZ_{q^{2}-1}$ of $G$, which also
contains a long-root parabolic subgroup
$P_{H} = q^{1+8} \cdot SL_{2}(q^{3}) \cdot \ZZ_{q-1}$ of $H$.

It is well known
that $V|_{Z}$ affords all the nontrivial linear characters $\lambda$
of the long-root subgroup $Z := Z(P')$ (which is elementary abelian of order $q^{2}$),
and the corresponding eigenspaces $V_{\lambda}$ are permuted regularly by the
torus $\ZZ_{q^{2}-1}$. Let $U = q^{2+16}$ denote the unipotent radical of $P$ and
consider any such $\lambda$.
Then $\IBRL(U)$ contains a unique representation (of degree $q^{8}$),
on which $Z$ acts via the character $\lambda$. Moreover, since
$P'/U \simeq SL_{2}(q^{6})$ has trivial Schur multiplier and is perfect, this
representation of $U$ extends to a unique representation of $P'$, which we denote
by $E_{\lambda}$. By Clifford theory, the $P'$-module $V_{\lambda}$ is isomorphic to
$E_{\lambda} \otimes A$ for some $A \in \IBRL(P'/U)$. Suppose that $A$ contains a
nontrivial composition factor, as a $SL_{2}(q^{6})$-module. Then
$\dim(A) \geq (q^{6}-1)/2$. It follows that
\begin{equation}\label{b1}
  \dim(V) \geq (q^{2}-1)q^{8}(q^{6}-1)/2.
\end{equation}
On the other hand, the irreducibility of $V|_{H}$ implies that
$$\dim(V) < \sqrt{|H|} < q^{14},$$
contradicting (\ref{b1}). Thus all composition factors of $A$ are trivial.
In particular, the $P'$-module $V_{\lambda}$ contains a simple submodule which is
isomorphic to $E_{\lambda}$.

Notice that we can embed $P_{H}$ in $P$ in such a way that $Z$ contains
$Z_{H} := Z(P'_{H})$ (a long-root subgroup in $H$, which is elementary abelian
of order $q$), and $U$ contains the unipotent radical $U_{H} = q^{1+8}$ of
$P_{H}$. Now choose $\lambda$ such that $Z_{H} \leq \Ker(\lambda)$. Then it is
easy to see that $E_{\lambda}|_{U_{H}}$ is just the regular representation, whence
the subspace $L$ of $U_{H}$-fixed points in it is one-dimensional, and, since
$U_{H} \lhd P'_{H}$, this subspace is acted on by $P'_{H}$. But
$P'_{H}/U_{H} \simeq SL_{2}(q^{3})$ is perfect, hence $P'_{H}$ acts trivially on
$L$.

We have shown that, for the given choice of $\lambda$, $P'_{H}$ has nonzero
fixed points in $V_{\lambda}$. Let $W$ be the subspace consisting of all
$P'_{H}$-fixed points in $V$. Then $P_{H}/P'_{H} \simeq \ZZ_{q-1}$ acts on $W$
and so $W$ contains a one-dimensional $P_{H}$-submodule $T$. By the Frobenius
reciprocity, $0 \neq \dim\Hom_{P_{H}}(T,V|_{P_{H}}) = \dim\Hom_{H}(\Ind^{H}_{P_{H}}(T),V|_{H})$.
But $V|_{H}$ is irreducible, hence it is a quotient of $\Ind^{H}_{P_{H}}(T)$. In particular,
\begin{equation}\label{b2}
  \dim(V) \leq (H:P_{H}) \cdot  \dim(T) = (q+1)(q^{8}+q^{4}+1).
\end{equation}
On the other hand, Theorem 4.1 of \cite{MMT} implies that
$$\dim(V) \geq q^{2}(q^{8}-q^{4}+1)-1,$$
contradicting (\ref{b2}).
\end{proof}

\section{Restriction to maximal parabolic subgroups}

\begin{lemma}\label{paraQ}
{\sl Let $Q$ denote the maximal parabolic subgroup of order
$q^{12}(q^{3}-1)(q^{2}-1)$ of $G = \tb D_{4}(q)$, and let $U := O_{p}(Q)$. Then

{\rm (i)} For any prime $r \neq p$, $O_{r}(Q) = 1$.

{\rm (ii)} Let $\varphi \in \IBRL(Q)$ be an irreducible Brauer
character of $Q$ whose kernel does not contain $U$. If $q$ is odd,
assume in addition that $\varphi$ is faithful. Then $\varphi$ lifts
to a complex character $\chi$ of $Q$. Moreover, $\chi$ is also
faithful if $\varphi$ is faithful.}
\end{lemma}

\begin{proof}
(i) Since $O_{r}(Q), U \lhd Q$ and $O_{r}(Q) \cap U = 1$, any element $g \in O_{r}(Q)$
is centralized by $U$, which has order $q^{11}$. Thus $q^{11}$ divides $|C_{Q}(g)|$.
Assuming $g \neq 1$, we see by \cite{H1} and \cite{H3} that $g$ is $Q$-conjugate to
the long-root element $u = x_{3\alpha+2\beta}(1)$. But then $g$ is a $p$-element,
a contradiction. Hence $O_{r}(Q) = 1$.

(ii) Let $\lambda$ be an irreducible constituent of $\varphi|_{U}$,
and let $I$ denote the inertia group of $\lambda$ in $Q$. By
Clifford theory, $\varphi = \Ind^{Q}_{I}(\psi)$ for some $\psi \in
\IBRL(I)$ whose restriction to $U$ contains $\lambda$. Since $p \neq
\ell$, we may view $\lambda$ as an ordinary character of $U$. By our
assumption, $\lambda \neq 1_{U}$. The structure of $I/U$ is
described in \cite{H1}, \cite{H3}. In particular, if $2|q$, then
$I/U$ is always solvable. On the other hand, if $q$ is odd, then
$I/U$ is solvable, except for one orbit, the kernel of any character
in which however contains a long-root element
$x_{3\alpha+2\beta}(1)$ (in the notation of \cite{H1}). Recall we are assuming
that $\varphi$ is faithful if $q$ is odd. It follows
that in either case $I/U$ is solvable, and so $I$ is solvable. By the Fong-Swan
Theorem, $\psi$ lifts to a complex character $\rho$ of $I$. Hence
$\varphi$ lifts to the complex character $\chi :=
\Ind^{Q}_{I}(\rho)$.

Now assume that $\varphi$ is faithful but $K := \Ker(\chi)$ is
non-trivial; in particular, $\ell \neq 0$. If $K$ is not an
$\ell$-group, then $K$ contains a non-trivial $\ell'$-element $g$.
Since $\varphi(g) = \chi(g) = \chi(1) = \varphi(1)$, we see that
$\varphi$ is not faithful, a contradiction. Hence $K$ is an
$\ell$-group, and so $O_{\ell}(Q) \neq 1$, contradicting (i).
\end{proof}

\begin{theor}\label{para}
{\sl Let $M$ be a maximal parabolic subgroup of $G = \tb D_{4}(q)$
and $\varphi \in \IBRL(G)$ be of degree $>1$. Then $\varphi|_{M}$ is
reducible.}
\end{theor}

\begin{proof}
First suppose that $M = P$, the long-root parabolic subgroup of $G$. Then
the statement follows from Theorem 1.6 of \cite{T}. So we may assume that
$M = Q$, the other maximal parabolic subgroup of $G$. Also assume the contrary:
$\varphi|_{Q}$ is irreducible.

We will consider two particular long-root elements $u =
x_{3\alpha+2\beta}(1)$ and $v = x_{\beta}(1)$ of $Q$, in the
notation of \cite{H1}, \cite{H2}, \cite{H3}. Clearly, they are
conjugate in $G$, so $\varphi(u) = \varphi(v)$. By Lemma
\ref{paraQ}, $\varphi|_{Q}$ lifts to a complex irreducible character
$\chi$ of $Q$ which is also faithful. Since $u$ and $v$ are
$\ell'$-elements, we have $\varphi(u) = \chi(u)$ and $\varphi(v) =
\chi(v)$. It follows that
\begin{equation}\label{uv}
  \chi(u) = \chi(v).
\end{equation}
Note that $Z := Z(O_{p}(Q)) = X_{3\alpha+2\beta}X_{3\alpha+3\beta}$
has order $q^{2}$, and consists of the $q^{2}-1$ $Q$-conjugates of $u$ and $1$. Thus
$Q$ acts transitively on $Z \setminus \{1\}$ and on $\Irr(Z) \setminus \{1_{Z}\}$.
Since $\Ker(\chi) = 1$, we conclude that $\chi(u) = -\chi(1)/(q^{2}-1)$.

First consider the case $q$ is odd. Then $u$, resp. $v$, belongs to the $Q$-conjugacy class
$c_{1,1}$, resp. $c_{1,2}$, in the notation of \cite{H1}. According to \cite{H1}, the
faithful character $\chi$ must be one of $\chi_{j}(k)$, $16 \leq j \leq 20$. If $j = 16$ or
$17$, then $\chi(v)$ is explicitly computed in \cite{H1}, and one sees that
(\ref{uv}) is violated.
Now suppose that $j = 18$ or $19$. Then $\chi(u) = -q^{3}(q^{3}-1)/2$. On the other hand,
according to Proposition 2.1 of the Appendix,
$\chi(v) = mq(q^{3}-1)$ with $m \geq -(q^{2}-1)/2$.
It follows that $\chi(v) > \chi(u)$, violating (\ref{uv}). Finally, suppose that
$j = 20$. Then $\chi(u) = -q^{3}(q^{3}-1)$. Meanwhile, by Proposition 2.1 of the Appendix,
$\chi(v) = mq(q^{3}-1)$ with $m \geq -(q^{2}-1)$.
It follows that $\chi(v) > \chi(u)$, again violating (\ref{uv}).

Next we consider the case $q$ is even. Then $u$, resp. $v$, belongs to the $Q$-conjugacy class
$c_{1,1}$, resp. $c_{1,7}$, in the notation of \cite{H3}. According to \cite{H3}, the
faithful character $\chi$ must be one of $\chi_{j}(k)$, $14 \leq j \leq 16$. If $j = 14$ or
$15$, then $\chi(u)$ is explicitly computed in \cite{H3}, and one sees that
(\ref{uv}) is violated. Finally, suppose that $j = 16$. Then $\chi(u) = -q^{3}(q^{3}-1)$.
On the other hand, by Proposition 1.1 of the Appendix,
$\chi(v) = mq(q^{3}-1)$ with $m \geq -(q^{2}-1)$.
It follows that $\chi(v) > \chi(u)$, again violating (\ref{uv}).
\end{proof}

{\bf Proof of Theorem \ref{main}.} Assume the contrary: $\Phi|_{H}$
is irreducible. Without loss we may assume that $\Phi$ is absolutely
irreducible and that $H$ is a maximal subgroup of $G$. Now we can
apply Theorem \ref{reduction} to $H$ to get four possibilities (i)
-- (iv) for $H$. None of them cannot however occur by Theorems
\ref{g2}, \ref{square}, and \ref{para}. \hfill $\Box$

\newcommand{\mn}{\mbox{-}}            
\newcommand{\mydelsep}{2.5pt}         
\newcommand{\mycolsep}{2.2pt}         

\newpage
\centerline{{\large {\bf
{Appendix}}}}
\vspace{0.3cm}
\centerline{{\large {\bf
{Faithful characters of a parabolic subgroup of $^3D_4(q)$}}}}

\bigskip
\centerline{by Frank Himstedt}

\bigskip
\bigskip
Let $q$ be a prime power, ${^3D}_4(q)$ Steinberg's triality group,
and $Q$ a maximal parabolic subgroup of order
$q^{12}(q^3-1)(q^2-1)$ of ${^3D}_4(q)$. A classification and
construction of all irreducible characters of $Q$ is described
in~\cite{Himstedt3D4Parab, Himstedt3D4Parab2},
and the values of almost all of these characters are given by
Tables~A.13 and A.14 in \cite{Himstedt3D4Parab,
Himstedt3D4Parab2}. This appendix deals with some values of the
faithful irreducible characters of $Q$ which are \emph{not} contained in
\cite{Himstedt3D4Parab, Himstedt3D4Parab2}. More specifically,
we are interested in the values of the faithful irreducible characters of $Q$
on the ``long root element''~$x_\beta(1)$ (for a definition of
$x_\beta(1)$, see \cite{GeckDissPub}, for example). In particular, we
prove bounds on these character values by considering scalar products
with restrictions of unipotent irreducible characters of ${^3D}_4(q)$.


Suppose that $q$ is even. We use the same notation as in
\cite{Himstedt3D4Parab2}. The faithful irreducible characters of~$Q$ are the
characters ${_Q\chi}_{14}(k)$, ${_Q\chi}_{15}(k)$, ${_Q\chi}_{16}(k)$ with $k$
as in \cite[Table~A.13]{Himstedt3D4Parab2}. The values of ${_Q\chi}_{14}(k)$
and ${_Q\chi}_{15}(k)$ are given in
\cite[Table~A.14]{Himstedt3D4Parab2}, in particular those on
$x_\beta(1)$. So we only deal with the values of ${_Q\chi}_{16}(k)$ on
$x_\beta(1)$.

\begin{propo} \label{prop:vals_even}
The values of $_Q\chi_{16}(k)$ on $x_\beta(1)$ satisfy:
\[
_Q\chi_{16}(k)(x_\beta(1)) \in \{ q (q^3-1) \cdot m \, | \, m \in \ZZ \text{
  with   } -(q^2-1) \le m \le q^2(q-1) \}.
\]
\end{propo}

\begin{proof}
We recall the definition of $_Q\chi_{16}(k)$ (see
\cite[p.~258 and p.~256]{Himstedt3D4Parab2}): number the elements of
the field $\FF_q$ in some way, say $\FF_q = \{x_1, x_2, \dots, x_q\}$ with
$x_1=0$. Then, $_Q\chi_{16}(k)$ is the character of $Q$
which is induced from the following linear character of the subgroup
$X_\beta X_{\alpha+\beta} X_{2\alpha+\beta} X_{3\alpha+\beta}
X_{3\alpha+2\beta}$:
\[
x_\beta(d_2) x_{\alpha+\beta}(d_3) x_{2\alpha+\beta}(d_4)
x_{3\alpha+\beta}(d_5) x_{3\alpha+2\beta}(d_6) \mapsto \phi'(x_k
\cdot d_2 + d_4+d_5)
\]
where $\phi'$ is a linear character of the additive group
of $\FF_{q^3}$ restricting nontrivially on~$\FF_q$.
Using the definition of induced characters and the relations
in Tables~2.2-2.4 in \cite{GeckDissPub} we see that the value of
$_Q\chi_{16}(k)$ on $x_\beta(1)$ is:
\begin{equation} \label{eq:faithval}
y_k := \sum_{i=1}^{q^3-1} \sum_{j=1}^{q-1} \sum_{s \in \FF_{q^3}}
\phi'(\tau^{(-q^2+q+1)i}s^{q+1} + x_k \pi^{2j-i} + \pi^{i-j} s^{q^2+q+1})
\end{equation}
where $\tau$ is a generator of the multiplicative group
$\FF_{q^3}^\times$ and $\pi := \tau^{q^2+q+1}$ a generator
of~$\FF_q^\times$.
Since $x_\beta(1)$
is an involution we know $y_k \in \ZZ$ and so $y_k \le q^3(q^3-1)(q-1)$.
Let $[\varepsilon_2]_Q$ be the restriction of the unipotent irreducible
character $[\varepsilon_2]$ of degree $q^7(q^4-q^2+1)$ of ${^3D}_4(q)$
(see~\cite{Spaltenstein}). Because we know the values of
$_Q\chi_{16}(k)$ on all conjugacy classes of $Q$ where
$[\varepsilon_2]_Q$ is nonzero, we can compute the scalar product
$(_Q\chi_{16}(k), [\varepsilon_2]_Q)_Q =
\frac{q^6-q^4-q^3+q+y_k}{q(q^3-1)}$ using {\sf CHEVIE}
\cite{CHEVIE}. Since this scalar product is a nonnegative rational
integer the statement about $_Q\chi_{16}(k)(x_\beta(1))$ follows.
\end{proof}

\noindent \textbf{Remark:} The value of $_Q\chi_{16}(1)$ on
$x_\beta(1)$ can be evaluated explicitly: using (\ref{eq:faithval}) and
$\sum_{i=1}^{q^3-1} \phi'(\tau^i) = \sum_{j=1}^{q-1} \phi'(\pi^j) =
-1$ one gets $_Q\chi_{16}(1)(x_\beta(1)) = q(q^3-1)$.

\medskip
Now suppose that $q$ is odd.
We use the same notation as in \cite{Himstedt3D4Parab}. The faithful
irreducible characters of~$Q$ are ${_Q\chi}_{16}(k)$,
${_Q\chi}_{17}(k)$, ${_Q\chi}_{18}(k)$, ${_Q\chi}_{19}(k)$,
${_Q\chi}_{20}(k)$ with $k$ as in
\cite[Table~A.13]{Himstedt3D4Parab}. The values of
${_Q\chi}_{16}(k)$, ${_Q\chi}_{17}(k)$ are given in
\cite[Table~A.14]{Himstedt3D4Parab}, in particular the values on
$x_\beta(1)$. So we only deal with the values of ${_Q\chi}_{18}(k)$,
${_Q\chi}_{19}(k)$ and ${_Q\chi}_{20}(k)$ on $x_\beta(1)$.

\begin{propo} \label{prop:vals_oddq}
The values ${_Q\chi}_{18}(k)(x_\beta(1))$ and
${_Q\chi}_{19}(k)(x_\beta(1))$ are elements of the set
$
\{ q(q^3-1) \cdot m \, | \, m \in \ZZ \text{   with   } -(q^2-1)/2 \le m \le
q^2(q-1)/2 \}
$
and
\[
{_Q\chi}_{20}(k)(x_\beta(1)) \in \{ q(q^3-1) \cdot m \, | \, m \in \ZZ \text{
  with   } -(q^2-1) \le m \le q^2(q-1) \}.
\]
\end{propo}
\begin{proof}
The proof is similar to that of Proposition~\ref{prop:vals_even}.
By construction (see p.~795 and 783 in \cite{Himstedt3D4Parab}), the
value of ${_Q\chi}_{18}(k)$ and ${_Q\chi}_{19}(k)$ on $x_\beta(1)$ is:
\[
\sum_{i=1}^{q^3-1} \sum_{j=1}^{(q-1)/2} \sum_{s \in \FF_{q^3}}
\phi'(-\tau^{(-q^2-q+1)i} \pi^j s - \pi^{i-j+k} s^{q^2+q+1}) =: y_k'.
\]
Because $x_\beta(1)$ is an involution we have $y_k' \in \ZZ$
and $y_k' \le q^3(q^3-1)(q-1)/2$. Using {\sf CHEVIE}, we compute the
scalar products $(_Q\chi_{18}(k), [\varepsilon_2]_Q)_Q = (_Q\chi_{19}(k),
[\varepsilon_2]_Q)_Q = \frac{q^6-q^4-q^3+q+2y_k'}{2q(q^3-1)}$.
The fact that this is a nonnegative integer implies the statement
about ${_Q\chi}_{18}(k)(x_\beta(1))$, ${_Q\chi}_{19}(k)(x_\beta(1))$.
By construction (see \cite[p.~795 and 784]{Himstedt3D4Parab}):
\[
{_Q\chi}_{20}(k)(x_\beta(1)) =
\sum_{i=1}^{q^3-1} \sum_{j=1}^{q-1} \sum_{s \in \FF_{q^3}}
\phi'(\pi^{2j-i+k} - \tau^{(-q^2-q+1)i} \pi^j s - \pi^{i-j}
s^{q^2+q+1}) =: y_k''
\]
for $1 \le k \le (q-1)/2$ and
\[
{_Q\chi}_{20}(k)(x_\beta(1)) =
\sum_{i=1}^{q^3-1} \sum_{j=1}^{q-1} \sum_{s \in \FF_{q^3}}
\phi'(\pi^{2j-i+k} - \tau^{(-q^2-q+1)i} \pi^j s - \pi^{i-j+1}
s^{q^2+q+1}) =: y_k''
\]
for $(q+1)/2 \le k \le q-1$. So $y_k'' \in \ZZ$ and
$y_k'' \le q^3(q^3-1)(q-1)$. The fact
$(_Q\chi_{20}(k), [\varepsilon_2]_Q)_Q =
\frac{q^6-q^4-q^3+q+y_k''}{q(q^3-1)} \in \ZZ_{\ge 0}$ completes the
proof.
\end{proof}


\end{document}